\documentclass[12pt]{amsart}
\usepackage{amssymb,amscd}
\usepackage[mathscr]{eucal}
\input epsf
\newtheorem{thm}{Theorem}[section]
\newtheorem{lemma}[thm]{Lemma}
\newtheorem{defin}[thm]{Definition}

\newtheorem{prop}[thm]{Proposition}

\newtheorem{cor}[thm]{Corollary}
\begin{document}

\newcommand{\C}{{\Bbb C}}
\newcommand{\HH}{{\Bbb H}}
\newcommand{\M}{{\mathcal M}}
\newcommand{\B}{{\mathcal B}}
\newcommand{\MB}{{\mathcal M}_{\mathcal B}}
\newcommand{\R}{{\Bbb R}}
\newcommand{\Z}{{\Bbb Z}}
\newcommand{\HB}{{\mathcal H}_\B}
\newcommand{\SU}{\o{SU}(2)}
\newcommand{\SLC}{\o{SL}(2,\C)}
\newcommand{\SLR}{\o{SL}(2,\R)}
\newcommand{\Hom}{\o{Hom}}
\newcommand{\Spin}{\o{Spin}(2)}
\newcommand{\Spinm}{\o{Spin}_-(2)}
\newcommand{\Pin}{\o{Pin}(2)}
\renewcommand{\O}{{\mathcal O}}
\renewcommand{\P}{{\Bbb P}}
\newcommand{\Q}{{\Bbb Q}}
\renewcommand{\o}{\operatorname}
\newcommand{\HomU}{\o{Hom}_\B(\pi_1(M), \SU)}
\newcommand{\HomR}{\o{Hom}_\B^+(\pi_1(M), \SLR)}
\newcommand{\HomC}{\o{Hom}_\B^+(\pi_1(M), \SLC)}
\newcommand{\HomG}{\o{Hom}_\B^+(\pi_1(M), G)}
\newcommand{\tr}{\o{tr}}
\title[Dynamics of mapping class group actions on Moduli Spaces]
{Dynamics of the Mapping Class Group on the Moduli of a Punctured Sphere
 with Rational Holonomy}
\author{Joseph P. Previte \and Eugene Z. Xia}
\address{ School of Science,
Penn State Erie, The Behrend College,
Erie, PA 16563}
\address{Department of Mathematics,
National Cheng Kung University, Tainan 701, Taiwan, R.O.C}
\email{jpp@vortex.bd.psu.edu ({\it Previte}), ezxia@ncku.edu.tw ({\it Xia})}
\date{\today}
\maketitle
\begin{abstract}
Let $M$ be a four-holed sphere and
$\Gamma$ the mapping class group of $M$ fixing the boundary $\partial M$.
The group $\Gamma$ acts on $\MB(\SLC) = \HomC/\SLC$ which
is the space of completely reducible $\SLC$-gauge equivalence classes of flat
$\SLC$-connections on $M$ with fixed holonomy $\B$ on
$\partial M$.
Let $\B \in (-2,2)^4$ and
$\MB$ be the compact component of the real points of
$\MB(\SLC)$.  These points
correspond to $\SU$-representations or
$\SLR$-representations.  The $\Gamma$-action preserves $\MB$ and 
we study the topological dynamics of the $\Gamma$-action on $\MB$ and show that
for a dense set of holonomy $\B \in (-2,2)^4$,
the $\Gamma$-orbits are dense
in $\MB$.  We also produce a class of representations $\rho \in \HomR$
such that
the $\Gamma$-orbit of $[\rho]$ is finite in the
compact component of $\MB(\SLR)$,
but $\rho(\pi_1(M))$ is dense
in $\SLR$.
\end{abstract}

\section{Introduction}
Let $M$ be a 4-holed sphere and
$\{\gamma_1, \gamma_2, \gamma_3, \gamma_4\} \subset \pi_1(M)$
be the elements in the fundamental group corresponding to
the four boundary components.  Denote by $\HomC$ the space of
completely reducible $\SLC$-representations.
Let $G$ be a semi-simple subgroup of $\SLC$.
Assign each $\gamma_i$ a conjugacy class in $G$:  For each
${\B}=(a,b,c,d) \in \C^4$, let
\begin{eqnarray*}
\HB(G) & = & \{\rho \in \HomG \subset \HomR : \\
& & (\tr(\rho(\gamma_1)),\tr(\rho(\gamma_2)),\tr(\rho(\gamma_3)),\tr(\rho(\gamma_4))) = \B \},
\end{eqnarray*}
where $\HomG \subset \HomC$ is the subspace of completely reducible 
$G$-representations.
We restrict ${\B}$ to be in $(-2,2)^4$.
The group $G$ acts on $\HB$
by conjugation.
\begin{defin}~\label{def:1.1}
The moduli space with fixed holonomy ${\B}$ is
$$
\MB(G) =  \HB(G)/G.
$$
Let $\MB$ be either $\MB(\SU)$ or $\MB(\SLR)^c$, the {\it compact
  component} of $\MB(\SLR)$, depending on $\B$ (see \cite{Be1}).
\end{defin}
Let $\o{Diff}(M, \partial M)$ be the group of diffeomorphisms of
$M$ fixing $\partial M$. The group $\Gamma$ acts on $\pi_1(M)$
naturally and the action induces an action: $\Gamma \times \HB (G)
\longrightarrow \HB (G)$ with $\gamma(\rho)(X) =
\rho(\gamma^{-1}(X))$.  This, in turn, gives an action $\Gamma
\times \MB \longrightarrow \MB$. 

The space $\MB$ possesses a natural
symplectic structure which gives rise
to a finite measure $\mu$ on $\MB$.  Goldman first showed that
$\Gamma$ acts ergodically on $\MB(\SU)$ (see \cite{Go1, Go2}).
The same technique of \cite{Go2} immediately gives ergodicity of the
$\Gamma$-action on $\MB$.

The space $\M_{\mathcal B}$ has a natural topology and one can study topological
dynamics of the $\Gamma$-action \cite{Pr1,Pr2,Pr3}.
%
In this paper, we show that for a dense subset of boundary holonomies,
the $\Gamma$-action is minimal:
\begin{thm}\label{thm:main}
Let $M$ be a four-holed sphere.  Suppose ${\mathcal B} = (a,b,c,d) \in
(-2,2)^4$ such that two of $(ab+cd), (ac+bd), (ad+bc)$ are in the set
$\Q \setminus
\{i \in \Z : -8 \le i \le 8\}$.
Then every $\Gamma$-orbit is dense in $\MB$.
\end{thm}

\begin{cor}
There exists a dense subset $D \subset [-2,2]^4$ such that ${\B}
\in D$ implies that the $\Gamma$-action is minimal in $\MB$.
\end{cor}

Suppose $G \subset \SU$ (resp. $G \subset \SLR$) is a closed proper
subgroup and $\rho \in
\HB(G)$.  Then the $\Gamma$-orbit of $[\rho]$ is in
$\HB(G)$.
Hence the $\Gamma$-action is not, in general, minimal.
However, if $M$ is a surface of positive genus and
if the image of $\rho \in \HB (\SU)$ 
is dense in $\SU$ then the $\Gamma$-orbit of
$[\rho]$ is dense in $\MB(\SU)$ \cite{Pr2}.
This is no longer true without the $g>0$ hypothesis as examples in
\cite{Pr3} illustrate. In this paper, we also produce similar examples
of $\rho \in \HB (\SLR)$ having dense images in $\SLR$ but with discrete
$\Gamma$-orbits in $\MB(\SLR)^c$.

We adopt the following notational conventions:
For a fixed $\rho \in \HB$ and $X \in \pi_1(M),$
we write $X$ for $\rho(X)$
when there is no ambiguity.
A small letter denotes the
trace of the matrix represented by the corresponding
capital letter.

\section{The $\Gamma$-action on moduli spaces}
We first review some results that appear
in \cite{Be1} and \cite{Go1}.
Suppose $M$ is a four-holed sphere.  Then the fundamental
group $\pi_1(M)$ admits a presentation
$$
\langle A, B, C, D : ABCD = I \rangle.
$$
Let $X = AB, Y = BC$ and  $Z = CA$.
For a fixed holonomy (trace) $\B \in (-2,2)^4$ on the four punctures
and with the above coordinates, the moduli space satisfies
\begin{eqnarray*}
x^2 + y^2 + z^2 + xyz & = &
(ab+cd)x + (ad + bc)y + (ac+bd)z \\
& & -(a^2 + b^2 + c^2 + d^2 + abcd -4).
\end{eqnarray*}
The $x$-level sets $X(x) \subset \MB$
satisfy:
$$
 \frac{2+x}{4}\biggl [(y+z) - \frac{(a+b)(d+c)}{2+x}\biggr ]^2 +
 \frac{2-x}{4}\biggl [(y-z) - \frac{(a-b)(d-c)}{2-x}\biggr ]^2 =
$$
$$
\frac{(x^2-abx+a^2+b^2-4)(x^2-cdx+c^2+d^2-4)}{4-x^2}.
$$
There are similar descriptions for the $y$- and $z$-level sets
$Y(y)$ and $Z(z)$ respectively (see \cite{Go1}).

Let
$$
I_{a,b}^- =  \frac { ab - \sqrt{(a^2-4)(b^2-4)}}2
$$
$$
I_{a,b}^+=
\frac {ab + \sqrt{(a^2-4)
(b^2-4)}}2
$$
(similarly for $I_{a,c}^+,I_{b,c}^-,$ etc.).

The moduli space $\MB$ is $\MB(\SU)$ if for some$(x,y,z)\in \MB$
$$
x^2-abx+a^2+b^2-4 < 0,
$$
$$
x^2-cdx+c^2+d^2-4 <0.
$$
These two inequalities imply that $I_{a,b}^- < I_{c,d}^+$ (or
$I_{c,d}^-< I_{a,b}^+$). Moreover, for $x \in S = (I_{a,b}^-,
I_{c,d}^+)$ (resp. in $(I_{c,d}^-, I_{a,b}^+)),$ the level set
$X(x)$ is an ellipse.

The moduli space $\MB$ is $\MB(\SLR)^c$ if for some $(x,y,z) \in
\MB,$
$$
x^2-abx+a^2+b^2-4 > 0,
$$
$$
x^2-cdx+c^2+d^2-4 > 0.
$$
In this case, $I_{a,b}^+ < I_{c,d}^-$ (or $I_{c,d}^+< I_{a,b}^-$).
Moreover, for $x \in S = (I_{a,b}^+, I_{c,d}^-)$ (resp. in
$(I_{c,d}^+, I_{a,b}^-)),$ the level set $X(x)$ is an ellipse.



%
%


For fixed $\B = (a,b,c,d)$, the
$x$-coordinates in $\MB$ take on all the values inside
$S$ (see \cite{Be1}).
In particular,
$$\MB = \bigcup_{x \in {\rm closure } {(S)}} X(x).$$
By symmetry, there exist similar constructions
for the $y$- and $z$-coordinates.

\subsection{The mapping class group action}

The mapping class group $\Gamma$ is generated by the maps $\tau_X$
and $\tau_Y$ induced by the Dehn twists in $X, Y \in \pi_1(M).$ In
local coordinates, the actions of $\tau_X, \tau_Y, \tau_Z$ are
$$\left[
\begin{array}{c}
y\\
z\\
\end{array}
\right]
\stackrel{\tau_X}{\longmapsto}
\left[
\begin{array}{c}
ad+bc-x(ac+bd-xy-z)-y\\
ac+bd- xy-z\\
\end{array}
\right],
$$
$$\left[
\begin{array}{c}
z\\
x\\
\end{array}
\right]
\stackrel{\tau_Y}{\longmapsto}
\left[
\begin{array}{c}
bd+ca-y(ba+cd-yz-x)-z\\
ba+cd-yz-x\\
\end{array}
\right],
$$
$$\left[
\begin{array}{c}
x\\
y\\
\end{array}
\right]
\stackrel{\tau_Z}{\longmapsto}
\left[
\begin{array}{c}
cd+ab- z(cb+ad-zx-y)-x\\
cb+ad- zx - y\\
\end{array}
\right].
$$
These three actions preserve the ellipses $X(x) \subset \MB,$ $Y(y)
\subset \MB$ and $Z(z) \subset \MB,$ respectively. After
coordinate transformations, these are rotations by angles
$2\cos^{-1}(x/2), 2\cos^{-1}(y/2)$ and $2\cos^{-1}(z/2)$,
respectively \cite{Go1}.

\section{The Irrational Rotations and Infinite Orbits}
The Dehn twist $\tau_Y$ acts on the (transformed)
subsets $Y(y)$ via a rotation of angle  $2\cos^{-1}(y/2)$.
Thus there is a filtration of the
$y$-coordinates that yields finite orbits under $\tau_Y$ \cite{Pr2}.
Let $Y_n \subset (-2, 2)$ such that $y \in Y_n$ if and only if
the $\tau_Y$-action on non-fixed points $(x,y,z) \in \MB$ is
periodic with period less than or equal to $n$.
This gives a filtration
$$
\{0\} = Y_2 \subset Y_3 = \{0, \pm 1\} \subset  Y_4 = \{0, \pm 1, \pm \sqrt 2\}$$

$$ \subset  Y_5 = \{0, \pm 1, \pm \sqrt 2, \pm \frac {1 \pm \sqrt 5}2 \} $$

$$\subset   Y_6 = \{0, \pm 1, \pm \sqrt 2, \pm \frac {1 \pm \sqrt 5}2 , \pm \sqrt 3\}
\subset ... \subset Y_n \subset ... $$
By symmetry, there are similar filtrations $X_n$ and $Z_n,$ with
$X_n=Y_n=Z_n$
as sets.

The global coordinates provide an embedding of $\MB$ in $\R^3.$ We consider
the box metric
$$
D((x_1,y_1,z_1), (x_2,y_2,z_2)) = max(|x_1 - x_2|, |y_1 - y_2|, |z_1 - z_2|),
$$
which generates the usual topology on $\MB$.

\begin{defin}
For $\epsilon > 0$, a set $U$ is
$\epsilon$-dense in $V$ if for each $p \in V \subset \MB$,
there exists a point $q \in U$ such that
$
0 < D(p, q) < \epsilon.
$
\end{defin}

\begin{lemma}\label{lem:path1}
For $\epsilon > 0$ there exists $N(\epsilon) > 0$
so that if $y \not \in Y_{N(\epsilon)}$,
then the $\tau_Y$-orbit of $(x,y,z)$
is $\epsilon$-dense in $Y(y)$ for any $(x,y,z)$ in any $\MB.$
\end{lemma}
\begin{proof}

Since the non-degenerate ellipses $Y(y)$
have uniformly bounded circumferences,
there exists $N(\epsilon) > 1$ such that
for any $y \not\in Y_{N(\epsilon)}$, the $\tau_Y$-orbit is $\epsilon$-dense
in $Y(y)$.
\end{proof}

Effectively, for any fixed $\epsilon$,
there is a finite number of $y$-coordinates (resp. $x$)
whose $\tau_Y$ (resp. $\tau_X$) actions are not
$\epsilon$-dense in $Y(y)$ (resp. $X(x)$).

\begin{thm}\label{thm:irrational}
If the $\Gamma$-orbit of $u_0 = (x_0, y_0, z_0)$ is infinite,
then it is dense in $\MB$.
\end{thm}
\begin{proof}

Let $u' = (x', y', z') \in \MB$ and
$\epsilon >0$.
By compactness, the $\Gamma$-orbit of $u_0$ has an accumulation
point $u_* = (x_*, y_*, z_*).$  Let $B_{\epsilon}(u_*)$ be the
ball of radius $\epsilon
> 0$ centered at $u_*$, with respect to the box metric.  Then either the set of
$x$-coordinates or the set of $y$-coordinates of the points in
$\Gamma (u_0) \cap B_{\epsilon}(u_*)$ are infinite. Since $\tau_X$
and $\tau_Y$ are rotations by angles $2\cos^{-1}(x/2)$ and
$2\cos^{-1}(y/2)$, by Lemma \ref{lem:path1},
$\Gamma (u_0) \cap B_{\epsilon}(u_*)$ contains an infinite subset  
whose points have their
$x-$coordinates and $y-$coordinates distinguished from each other.

Recall that
$X(x)$ (resp., $Y(y)$) is the $x-$cross section (resp.
$y-$cross section) of $\MB,$ which is topologically a circle
or a point.
Let
$$X_\epsilon (x)= \bigcup_{|r|<\epsilon} X(x+r).$$

By compactness, there is a finite chain of sets
$S_0= X_\epsilon (x_*)$, $S_1 = Y_\epsilon (y_1),$
$S_2 = X_\epsilon (x_2)$, ... , $S_n = Y_\epsilon (y')$
with $S_{i}\cap S_{i+1} \neq \emptyset$.

As $\Gamma (u_0)$ contains an infinite number of points in
$S_0$ with distinct $x-$coordinates, there is (by using $\tau_X$) an
infinite number of points in $\Gamma (u_0) \cap S_0 \cap S_1$ having distinct
$y-$coordinates. Continuing in this fashion,
one generates an infinite number of points inside
$\Gamma (u_0) \cap S_{i}\cap S_{i+1}$ which leads to an infinite number of points
inside $\Gamma (u_0) \cap B_{\epsilon}(u').$
\end{proof}

\section{Minimality}

By Theorem~\ref{thm:irrational}, the problem amounts to ensuring that
some Dehn twist corresponds to an
irrational rotation along a circle in $\MB$.
For a given representation $\rho$, we
use subscript notation to denote the
actions on the coordinates $x,y,z$ of $[\rho]$.
For instance,
$x_y$ is short for the $x$-coordinate of $\tau_Y([\rho])$.
Since $x_y  = ab+cd-yz-x$,
for the orbit to be finite,
$2\cos^{-1}(\frac{x_z}2)$ and $2\cos^{-1}(\frac{z_y}2)$
must also be rational multiples of $\pi.$
Hence
$$
\cos (\theta_{x_y})+ 2 \cos (\theta_z) \cos (\theta_{y}) +
\cos(\theta_x) = \frac{ab+cd} 2,
$$
or equivalently,
\begin{equation}\label{eqn:cos}
\cos (\theta_{x_y})+ \cos (\theta_{z}+\theta_{y}) +\cos
(\theta_{z}-\theta_{y}) + \cos(\theta_x) = \frac{ab+cd}2,
\end{equation}
where all angles are rational multiples of $\pi$,
$0\le \theta_z+\theta_{y} \le 2\pi$,
and $-\pi \le \theta_z-\theta_{y} \le \pi.$
We also obtain two more equations from the $\tau_X$ and $\tau_Z$ actions.

Equation (\ref{eqn:cos})
is a trigonometric Diophantine equation,
the solutions to which
are few as shown by Conway and Jones:
\begin{thm} [Conway, Jones] \cite{Co} \label{thm:3}
Suppose that we have at most four distinct rational multiples of $\pi$
lying strictly between $0$ and $\pi/2$ for which some linear
combination of their cosines is rational, but no proper subset
has this property. That is,
$$
A \cos( a) + B \cos(b) + C \cos( c) + D \cos(d) = E,
$$
for $A, B, C, D, E$ rational and $a,b,c,d \in (0, \pi/2)$
rational multiples of $\pi$.
Then the appropriate linear combination is proportional
to one from the following list:
\begin{eqnarray*}
\cos(\pi/3) & = & 1/2\\
\cos(t+\pi/3)+\cos(\pi/3-t)-\cos(t) & = & 0  \ \ (0< t <\pi/6)\\
\cos(\pi/5) - \cos(2\pi/5) & = & 1/2\\
\cos(\pi/7) - \cos(2\pi/7) + \cos(3\pi/7) & = & 1/2\\
\cos(\pi/5) - \cos(\pi/15) + \cos(4\pi/15) & = & 1/2\\
-\cos(2\pi/5)+\cos(2\pi/15)-\cos(7\pi/15) & = &  1/2\\
\cos(\pi/7) + \cos(3\pi/7) - \cos(\pi/21) +\cos(8\pi/21) & = & 1/2\\
\cos(\pi/7) - \cos(2\pi/7) + \cos(2\pi/21) -\cos(5\pi/21) & = & 1/2\\
-\cos(2\pi/7) + \cos(3\pi/7) + \cos(4\pi/21) +\cos(10\pi/21) & = & 1/2\\
-\cos(\pi/15) + \cos(2\pi/15) + \cos(4\pi/15) -\cos(7\pi/15) & = & 1/2.\\
\end{eqnarray*}
\end{thm}

The angles in equation (\ref{eqn:cos}) are not necessarily
in $(0,\pi/2)$.
By applying the identities $\cos(\pi/2-t)= - \cos(\pi/2+t)$
and $\cos(\pi-t)=\cos(\pi+t)$, we derive from
equation (\ref{eqn:cos}) a new four-term cosine equation whose
arguments are in $[0,\pi/2]$.
That is, by a possible change of sign,
each term in equation (\ref{eqn:cos}) may be rewritten with angles
in $[0,\pi/2].$
Notice that if the resulting equation has non-distinct angles or an
angle being $0$ or $\frac{\pi}{2}$, then we will obtain a
rational trigonometric Diophantine equation of shorter length.

\begin{prop}\label{prop:weak}
Suppose two of
${ab+cd},ac+bd,ad+bc\in \Q \setminus \{i \in \Z : -8 \le i \le 8\}$.  Then
Equation~(\ref{eqn:cos}) has no rational (in the sense of
Conway-Jones) solution.
\end{prop}
\begin{proof}
Let $(x,y,z) \in {\mathcal M}_{\mathcal B}$.  Since $\B \in (-2,2)^4$,
$\MB$ is a smooth 2-sphere.  Hence two of the level sets $X(x),
Y(y), Z(z)$ are circles.  We assume that
$Y(y)$ is a circle and $ac+bd \in \Q \setminus \{i \in \Z: -8
\le i \le 8\}$.  Then a finite $\tau_Y$-orbit on
$Y(y)$ must yield a rational solution (in the sense of Conway-Jones)
to Equation~(\ref{eqn:cos}).
However, with the given constraint on $\frac{ab+cd}2$,
Theorem~\ref{thm:3}
implies that
Equation~(\ref{eqn:cos}) has no solution.  The other two cases are similar.
\end{proof}
Theorem~\ref{thm:main} follows immediately from Proposition~\ref{prop:weak}.

\section{Exceptional discrete orbits}
In this section, we construct representations 
$\rho \in \HB (\SLR)$ such that
$\rho(\pi_1(M))$ is dense in $\SLR$, but the $\Gamma$-orbit of
$[\rho]$ is discrete in $\MB(\SLR)^c$.  The construction closely
parallels the construction in \cite{Pr3}.

Let $F$ be the set of
$\B = (a,a,c,-c) \in (-2,2)^4$ satisfying the following conditions:
\begin{enumerate}
\item $a^2 + c^2 > 4$,
\item $\frac{\cos^{-1}(a)}{\pi}$ or $\frac{\cos^{-1}(c)}{\pi}$ is irrational.
\end{enumerate}
Consider the space $\MB$ with $\B \in F$.
The orbit
$$\O = \{(a^2-2,0,0), (2-c^2,0,0)\} \subset \MB$$
is $\Gamma$-invariant.  Let $\rho \in \HB(\SLR)$ with $[\rho]$ having the
coordinate $(a^2-2,0,0)$ or $(2-c^2,0,0)$.
Then condition (1) guarantees that $\rho(\pi_1(M))$
is not abelian since it has a non-trivial $\Gamma$-action.
Condition (2) implies that
$\rho(\pi_1(M))$ is dense in $S^1 \subset \SLR$.  Since $-2 < a,c <
2$, $\rho(\pi_1(M))$ is not contained in the group of affine
transformations of the real line.  Hence $\rho(\pi_1(M))$ is dense in
$\SLR$.

For a concrete example of one such case, let
$\B = (1, 1, \frac{7}{4}, -\frac{7}{4})$.  The special
orbit $\O$ consists of the two points that are intersections of the $x$-axis
with $\MB$, i.e.
$\O=\{ (-1,0,0), (-\frac{17}{16}, 0, 0)\}.$
Below is a representation in the conjugacy class
$(-1,0,0) \in \O \subset \MB$:
$$A = B =
\left[
\begin{array}{cc}
\frac{4}{5} & - \frac{3}{5} \\
& \\
\frac{7}{5} & \frac{1}{5}\\
\end{array}
\right]  \ \ \mbox{and} \ \ \
C =
\left[
\begin{array}{cc}
1 & -\frac{1}{4} \\
& \\
1 &  \frac{3}{4} \\
\end{array}
\right]. \ \ $$


\end{document}